\documentstyle[amscd,amssymb,verbatim,epsf]{amsart}

\begin{document}
\theoremstyle{plain}
\newtheorem{Thm}{Theorem}
\newtheorem{Cor}{Corollary}
\newtheorem{Con}{Conjecture}
\newtheorem{Main}{Main Theorem}
\newtheorem{Lem}{Lemma}
\newtheorem{Prop}{Proposition}

\theoremstyle{definition}
\newtheorem{Def}{Definition}
\newtheorem{Note}{Note}
\newtheorem{Ex}{Example}

\theoremstyle{remark}
\newtheorem{notation}{Notation}
\renewcommand{\thenotation}{}

\errorcontextlines=0
\numberwithin{equation}{section}
\renewcommand{\rm}{\normalshape}%

\title%
   {Reflection of a wave off a surface}
\author{Brendan Guilfoyle}
\address{Brendan Guilfoyle\\
          Department of Mathematics and Computing \\
          Institute of Technology, Tralee \\
          Clash \\
          Tralee  \\
          Co. Kerry \\
          Ireland.}
\email{brendan.guilfoyle@@ittralee.ie}

\author{Wilhelm Klingenberg}
\address{Wilhelm Klingenberg\\
 Department of Mathematical Sciences\\
 University of Durham\\
 Durham DH1 3LE\\
 United Kingdom.}
\email{wilhelm.klingenberg@@durham.ac.uk }

\date{July 22, 2002}

\begin{abstract}
Recent advances in twistor theory are applied to geometric optics in
${\Bbb{R}}^3$. The general formulae for reflection of a wavefront in
a surface are derived and in three special cases explicit descriptions
are provided: when the
reflecting surface is a plane, when the incoming wave is a plane
and when the incoming wave is spherical. In each case particular
examples are computed exactly and the results plotted to illustrate
the outgoing wavefront.

\end{abstract}

\maketitle

\section{Introduction}

In geometric optics, Huygens' principle allows one to describe the
propagation of light from two alternative perspectives: one can trace
the rays or one can trace the wavefronts. The drawback with the former
description is that following a finite number of rays may not describe
all of the phenomenon under study, while the formation of caustics in
wavefronts can cause difficulties in the latter description \cite{arn}.

In this paper we utilise recent work \cite{gak1} in twistor theory to
go back and forth between these two perspectives and thus describe the most
elementary of optical phenomena: reflection. In particular, we
consider the following situation: an incoming wave of light is
reflected on a surface in ${\Bbb{R}}^3$.  Given the shape
of the incoming wavefront and the reflecting surface, can one describe
the reflected wavefront? Throughout, we assume that the medium is
homogenous and that the speed of propagation is unity. 

The technique we employ in answering this question comes from the
minitwistor correspondence which identifies the space of
oriented affine lines in ${\Bbb{R}}^3$ with the tangent bundle to the
2-sphere. This has a long history and has been used in various
contexts. In particular,  it has been used in the construction of
minimal surfaces \cite{weier}, solutions to the wave equation
\cite{whitt} and the monopole equation \cite{hitch}. 

The general context of this work is within the study of line
congruences, that is, 2-parameter families of oriented lines in
${\Bbb{R}}^3$ \cite{hlav} \cite{hosch} \cite{pott}. A line congruence
is integrable if it is orthogonal to a family of surfaces in
${\Bbb{R}}^3$.  At the outset it is not clear that the reflection of
an integrable congruence in an arbitrary surface is itself
integrable. However, in Theorem \ref{t:int} that is precisely what we
prove. This is the celebrated Theorem of Malus, independently proven
by Hamilton in 1827.

In our formalism, going from an integrable line
congruence to the orthogonal surface is equivalent to inverting a
Cauchy-Riemann operator. The
reflection problem can then be broken down into four steps. First we
transfer from the incoming wavefront to its associated ray system and
we do the same for the reflecting surface. The reflection law
relates the outgoing rays to these two line congruences at the point of
contact, and, finally, we transform the outgoing rays to the
associated wavefront.

We present the general formulae relating all of these stages for an
arbitrary incoming wavefront and reflecting surface (Theorem
\ref{t:main}). We go on to give explicit descriptions for three special
cases: when the reflecting surface is a plane (Proposition
\ref{t:planesurf}), when the incoming wave is a plane (Proposition
\ref{t:planepar}) and when the incoming wave is spherical (Proposition
\ref{t:sphwave}). In each case, particular examples are computed
exactly and the results plotted to illustrate the outgoing wavefront.

The next section contains the background details of the twistor
construction we use - further details can be found in
\cite{gak2}.  In section 3 we deduce the required reflection law and
establish the fact that the outgoing ray system is integrable if and
only if the incoming rays are integrable. 

In section 4 we turn to the
simplest type of reflection, namely, reflection of an arbitary
wavefront in a plane. Section 5
deals with the general formulae for reflection of a plane wavefront in
an arbitrary surface. As an illustration of our technique we determine
the reflection of a plane wave in a sphere and a torus. In the final
section we give the general formulae for reflection of a spherical
wavefront in an arbitrary surface, and illustrate the technique in the
case where the surface is the unit sphere.

\section{The Twistor Correspondence}   

The key to our approach is the following geometric construction. Let
($x^1,x^2,x^3$) be Euclidean coordinates on
${\Bbb{R}}^3={\Bbb{C}}\oplus{\Bbb{R}}$ and set $z=x^1+ix^2$, $t=x^3$. 

Consider an oriented line in ${\Bbb{R}}^3$ passing through ($z$, $t$) with
direction given by $\xi\in S^2\subset{\Bbb{R}}^3$, where $\xi$ is
obtained by stereographic projection from the south pole of the unit
sphere about the origin onto the
plane through the equator.

Then the minimal distance vector from the line to the origin is given by
\cite{gak1}
\[
\eta\frac{\partial}{\partial \xi}+\overline{\eta}\frac{\partial}{\partial \overline{\xi}}\in T_\xi S^2,
\]
where 
\[
\eta=\frac{1}{2}(z-2t\xi-\overline{z}\xi^2),
\]
and the point ($z$, $t$) is a distance 
\[
r=\frac{\overline{\xi}z+\xi\overline{z}+(1-\xi\overline{\xi})t}{1+\xi\overline{\xi}}.
\]
from the point on the line closest to the origin.

Conversely, given $\eta$ and $r$ and the direction of the line $\xi$,
the point ($z$, $t$) in ${\Bbb{R}}^3$ can be found by
\begin{equation}\label{e:coord1}
z=\frac{2(\eta-\overline{\eta}\xi^2)+2\xi(1+\xi\overline{\xi})r}{(1+\xi\overline{\xi})^2},
\end{equation}
\begin{equation}\label{e:coord2}
t=\frac{-2(\eta\overline{\xi}+\overline{\eta}\xi)+(1-\xi^2\overline{\xi}^2)r}{(1+\xi\overline{\xi})^2}.
\end{equation}

Now an oriented line in ${\Bbb{R}}^3$ is uniquely determined by $\xi$ and
$\eta$ as above. This is the minitwistor correspondence \cite{hitch},
where the space of oriented lines in ${\Bbb{R}}^3$ is identified with
the tangent bundle to the unit sphere about the origin.

\begin{Def}
A {\it line congruence} is a 2-parameter family of oriented
lines. 
\end{Def}
 \begin{Def}
A line congruence is {\it integrable} if it is orthogonal to a family of
surfaces in ${\Bbb{R}}^3$.
\end{Def}

Suppose the line congruence is given by
$\nu\rightarrow(\xi(\nu,\bar{\nu}),\eta(\nu,\bar{\nu}))$, for
$\nu\in{\Bbb{C}}$. 

\begin{Thm}\cite{gak2}
A line congruence is integrable iff
\begin{equation}\label{e:intcond}
\partial\left(\frac{\eta}{(1+\xi\bar{\xi})^2}\bar{\partial}\bar{\xi}+
  \frac{\bar{\eta}}{(1+\xi\bar{\xi})^2}\bar{\partial}\xi \right)=
\bar{\partial}\left(\frac{\bar{\eta}}{(1+\xi\bar{\xi})^2}\partial\xi+
  \frac{\eta}{(1+\xi\bar{\xi})^2}\partial\bar{\xi} \right),
\end{equation}
where $\partial$ is differentiation with respect to $\nu$.

This equation is the integrability condition for the existence of a
real solution $r$ to the following equation 
\begin{equation}\label{e:pot}
\bar{\partial}r=\frac{2(\eta\bar{\partial}\bar{\xi}+\bar{\eta}\bar{\partial}\xi)}{(1+\xi\bar{\xi})^2},
\end{equation}
where $r$ is now considered as a function of $\nu$ and $\bar{\nu}$.
\end{Thm}

Thus, by the above theorem, given an integrable line congruence, we
can find a local description for the orthogonal surfaces by inverting
the $\bar{\partial}$ operator. The real constant of integration gives
the affine parameter along the lines normal to these surfaces. An
explicit description of the
surfaces in ${\Bbb{R}}^3$ can then be obtained by inserting $\xi$, $\eta$ and $r$, as
functions of $\nu$ and $\bar{\nu}$ into (\ref{e:coord1}) and
(\ref{e:coord2}).

Away from flat points, a surface can be parameterized by its normal
direction \cite{gak2}, i.e. $\nu=\xi$ and $\eta=F(\xi,\bar{\xi})$ for
some complex function $F$. In this case we refer to $F$ as the {\it
twistor function} of the surface, the integrability condition
(\ref{e:intcond}) reduces to the simpler
\[
\partial\left(\frac{F}{(1+\xi\bar{\xi})^2}\right)=
\bar{\partial}\left(\frac{\bar{F}}{(1+\xi\bar{\xi})^2} \right),
\]
and the function $r$ satisfies
\[
\bar{\partial}r=\frac{2F}{(1+\xi\bar{\xi})^2}.
\]
We will refer to $r$ as the {\it potential function} for the surface.

At a number of points in this paper we use Euclidean motions to
simplify the equations for reflection. A translation which takes the
origin to ($z$, $t$) is the quadratic holomorphic transformation
\[
\xi\rightarrow\xi, \qquad \qquad \eta\rightarrow\eta
    +\frac{1}{2}(z-2t\xi-\bar{z}\xi^2),
\]
while a rotation about the origin is given by the fractional linear
transformation:
\[
\xi\rightarrow\frac{\alpha\xi-\bar{\beta}}{\beta\xi+\bar{\alpha}},
\qquad \qquad \eta\rightarrow\frac{\eta}{(\beta\xi+\bar{\alpha})^2},
\]
where $\alpha, \beta\in{\Bbb{C}}$ satisfy
$\alpha\bar{\alpha}+\beta\bar{\beta}=1$.

\section{Reflection in a Surface}

Consider a ray with direction $\xi_1$ and perpendicular distance
vector $\eta_1$ striking a surface $S$ at the point ($z_0$,
$t_0$), where the normal direction is $\xi_0$ and  the perpendicular distance
vector  is $\eta_0$. Let $\xi_2$ be the direction of the outgoing wave and
$\eta_2$ be the perpendicular distance vector (see Figure 1). 

\setlength{\epsfxsize}{3in}
\begin{center}
   \mbox{\epsfbox{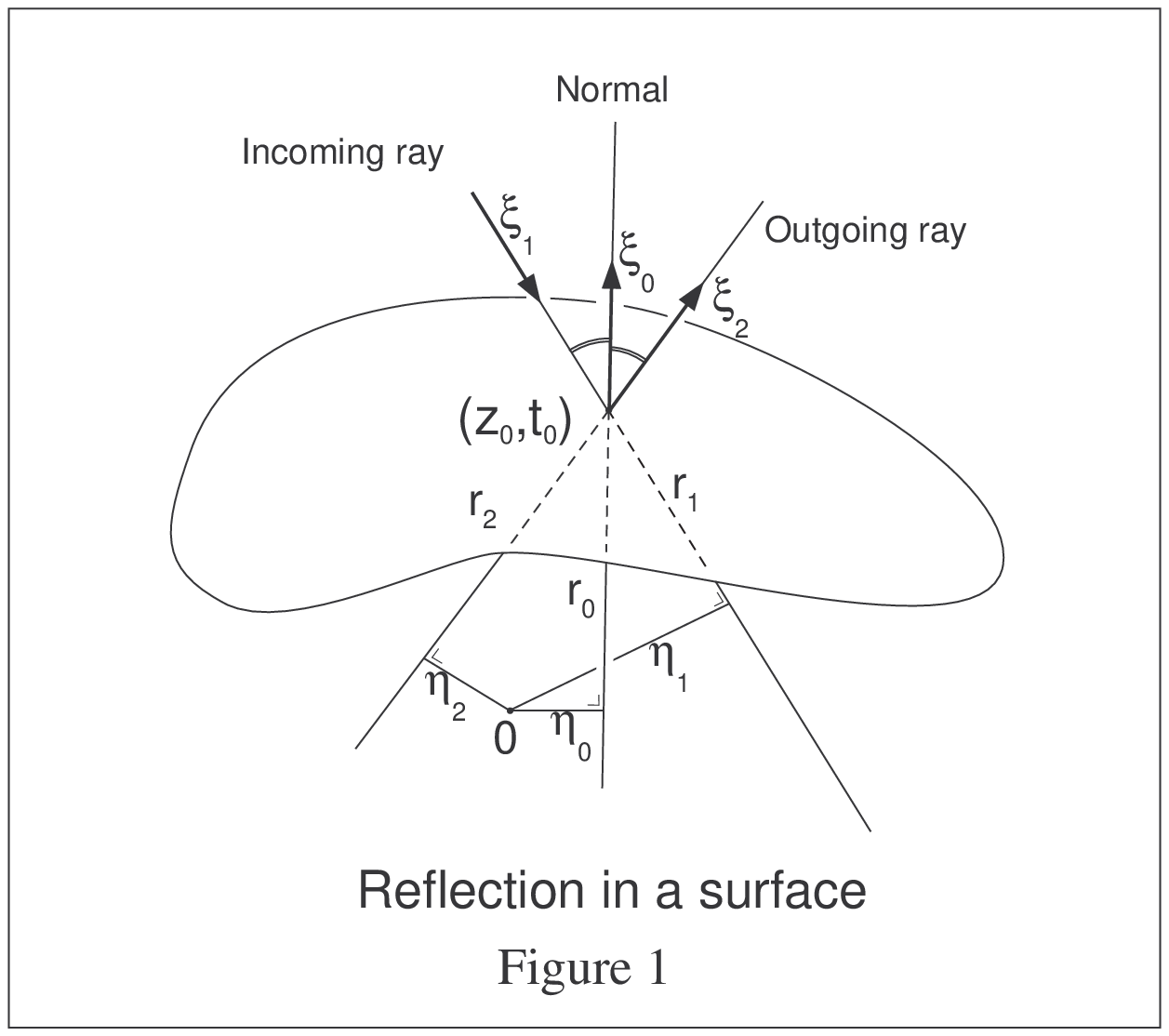}}
\end{center}

For the surface $S$, suppose that the point of reflection is a
distance $r_0$ from the minimal distance point on the normal line to the
origin, so that, by (\ref{e:coord1}) and (\ref{e:coord2})

\[
z_0=\frac{2(\eta_0-\overline{\eta}_0\xi_0^2)+2\xi_0(1+\xi_0\overline{\xi}_0)r_0}{(1+\xi_0\overline{\xi}_0)^2}, \qquad\qquad
t_0=\frac{-2(\eta_0\overline{\xi}_0+\overline{\eta}_0\xi_0)+(1-\xi_0^2\overline{\xi}_0^2)r_0}{(1+\xi_0\overline{\xi}_0)^2}.
\]

On the other hand, since both the incoming and outgoing rays contain
($z_0$, $t_0$),
\[
\eta_i=\frac{1}{2}(z_0-2t_0\xi_i-\overline{z}_0\xi_i^2),
\]
for $i=1,2$.

Combining the previous three equations we have, after some
rearrangement,
\begin{equation}\label{e:perpdists}
\eta_i=\frac{(1+\bar{\xi}_0\xi_i)^2\eta_0
-(\xi_0-\xi_i)^2\bar{\eta}_0
+(\xi_0-\xi_i)(1+\bar{\xi}_0\xi_i)(1+\xi_0\bar{\xi}_0)r_0
}{(1+\xi_0\bar{\xi}_0)^2}.
\end{equation}
We turn now to the law of reflection.

\begin{Prop}
If $\xi_1$ is the direction of the incoming ray, $\xi_2$ the outgoing
ray direction and  $\xi_0$ the normal direction at the point of
reflection, then

\begin{equation}\label{e:reflaw}
\xi_2=\frac{2\xi_0\bar{\xi}_1+1-\xi_0\bar{\xi}_0}
           {(1-\xi_0\bar{\xi}_0)\bar{\xi}_1-2\bar{\xi}_0}.
\end{equation}

\end{Prop}

\begin{pf}
A rotation about any point $\xi_0$ on ${\Bbb{P}}^1$ is described by a unitary
fractional linear transformation:
\[
\mbox{Rot}_{\xi_0}(\xi)=\frac{\alpha\xi-\bar{\beta}}{\beta\xi+\bar{\alpha}},
\]
where $\alpha, \beta\in{\Bbb{C}}$ satisfy
$\alpha\bar{\alpha}+\beta\bar{\beta}=1$. The inverse rotation is given
by
\begin{equation}\label{e:invrot}
\mbox{Rot}_{\xi_0}^{-1}(\xi)=\frac{\bar{\alpha}\xi+\bar{\beta}}{-\beta\xi+\alpha}.
\end{equation}
Let $S_{\xi_0}$ be the rotation
\[
S_{\xi_0}(\xi)=\frac{\alpha\xi-\bar{\beta}}{\beta\xi+\bar{\alpha}},
\]
where 
\[
\alpha=\frac{1}{\sqrt{1+\xi_0\bar{\xi}_0}} \qquad\qquad 
\beta=\frac{\bar{\xi}_0}{\sqrt{1+\xi_0\bar{\xi}_0}}.
\]
Clearly $S_{\xi_0}$ is a rotation that takes $\xi_0$ to zero (the
North pole). Thus, if we denote reflection through $\xi_0$ by
$R_{\xi_0}$, then $R_{\xi_0}=S_{\xi_0}^{-1}\circ R_{0}\circ S_{\xi_0}$.

Now $R_{0}(\xi)=-\xi$, so using (\ref{e:invrot}) and the definition of
$S_{\xi_0}$ we find that
\[
R_{\xi_0}(\xi)=\frac{(\xi_0\bar{\xi}_0-1)\xi+2\xi_0}
          {2\bar{\xi}_0\xi+1-\xi_0\bar{\xi}_0}.
\]
Alternatively, we can define $R_{\xi_0}$ as above and check that it is
a rotation (obvious), and satisfies $R_{\xi_0}(\xi_0)=\xi_0$ and
$R_{\xi_0}\circ R_{\xi_0}(\xi)=\xi$. Thus it is a 
rotation through 180$^0$ about $\xi_0$ i.e. reflection in 
$\xi_0$.

Finally, the law of reflection says that the direction $\xi_2$ of
the outgoing ray is obtained by reflecting the antipodal direction of
the incoming ray $\xi_1$ through the normal direction $\xi_0$ at the
point of reflection. 

To complete the proposition, we note that the antipodal map on
${\Bbb{P}}^1$ is $\xi\rightarrow-\bar{\xi}^{-1}$, so that
\[
\xi_2=R_{\xi_0}\left(-\frac{1}{\bar{\xi}}\right)
=\frac{2\xi_0\bar{\xi}_1+1-\xi_0\bar{\xi}_0}
           {(1-\xi_0\bar{\xi}_0)\bar{\xi}_1-2\bar{\xi}_0},
\]
as claimed.
\end{pf}

The equations governing reflection are:

\begin{Thm}\label{t:main}
Consider a wavefront given by
$\nu_1\rightarrow(\xi_1(\nu_1,\bar{\nu}_1),\eta_1(\nu_1,\bar{\nu}_1))$
reflecting off a surface given by
$\nu_0\rightarrow(\xi_0(\nu_0,\bar{\nu}_0),\eta_0(\nu_0,\bar{\nu}_0))$
and $r_0(\nu_0,\bar{\nu}_0)$. Then the reflected wavefront is
determined by

\begin{equation}\label{e:reflawa}
\xi_2=\frac{2\xi_0\bar{\xi}_1+1-\xi_0\bar{\xi}_0}
           {(1-\xi_0\bar{\xi}_0)\bar{\xi}_1-2\bar{\xi}_0},
\end{equation}
\begin{equation}\label{e:key1}
\eta_1=\frac{(1+\bar{\xi}_0\xi_1)^2\eta_0
-(\xi_0-\xi_1)^2\bar{\eta}_0
+(\xi_0-\xi_1)(1+\bar{\xi}_0\xi_1)(1+\xi_0\bar{\xi}_0)r_0,
}{(1+\xi_0\bar{\xi}_0)^2}
\end{equation}
\begin{equation}\label{e:key2}
\eta_2=\frac{(\bar{\xi}_0-\bar{\xi}_1)^2\eta_0
-(1+\xi_0\bar{\xi}_1)^2\bar{\eta}_0
+(\bar{\xi}_0-\bar{\xi}_1)(1+\xi_0\bar{\xi}_1)(1+\xi_0\bar{\xi}_0)r_0
}{((1-\xi_0\bar{\xi}_0)\bar{\xi}_1-2\bar{\xi}_0)^2}.
\end{equation}
Equation (\ref{e:key1})
determines the intersection of the incident rays with the surface,
while (\ref{e:reflawa}) and (\ref{e:key2}) determine the direction and
perpendicular distance from the origin of the reflected rays.
\end{Thm}
\begin{pf}
These come from combining the reflection law (\ref{e:reflaw}) with
(\ref{e:perpdists}), after some rearrangement for $\eta_2$.
\end{pf}

In general, the reflected line congruence can be parameterized by
$\nu_1$, the parameter for the incoming wavefront.

\begin{Thm}\label{t:int}
The surface, incident and reflected congruences satisfy the following
relationship:
\begin{equation}\label{e:twist}
\frac{\eta_2\bar{\partial}_1\bar{\xi}_2+\bar{\eta}_2\bar{\partial}_1\xi_2}{(1+\xi_2\bar{\xi}_2)^2}=
\frac{\eta_1\bar{\partial}_1\bar{\xi}_1+\bar{\eta}_1\bar{\partial}_1\xi_1}{(1+\xi_1\bar{\xi}_1)^2}+
\bar{\partial}_1\left(\frac{|\xi_0-\xi_1|^2-|1+\xi_0\bar{\xi}_1|^2}
    {(1+\xi_0\bar{\xi}_0)(1+\xi_1\bar{\xi}_1)}\;r_0
\right),
\end{equation}
where $\partial_1$ is differentiation with respect to $\nu_1$, the
parameter for the incoming wave.
\end{Thm}
\begin{pf}
We start by differentiating the reflection law (\ref{e:reflawa}):
\[
\partial_1\xi_2=\frac{-(1+\xi_0\bar{\xi}_0)^2\partial_1\bar{\xi}_1+2(\bar{\xi}_0- \bar{\xi}_1)^2\partial_1\xi_0+2(1+\xi_0\bar{\xi}_1)^2\partial_1\bar{\xi}_0}{[(1-\xi_0\bar{\xi}_0)\bar{\xi}_1-2\bar{\xi}_0]^2},
\]
\[
\bar{\partial}_1\xi_2=\frac{-(1+\xi_0\bar{\xi}_0)^2\bar{\partial}_1\bar{\xi}_1+2(\bar{\xi}_0-\bar{\xi}_1)^2\bar{\partial}_1\xi_0+2(1+\xi_0\bar{\xi}_1)^2\bar{\partial}_1\bar{\xi}_0}{[(1-\xi_0\bar{\xi}_0)\bar{\xi}_1-2\bar{\xi}_0]^2}.
\]
The reflection law (\ref{e:reflawa}) also implies
\[
1+\xi_2\bar{\xi}_2=\frac{(1+\xi_0\bar{\xi}_0)^2(1+\xi_1\bar{\xi}_1)}{[(1-\xi_0\bar{\xi}_0)\bar{\xi}_1-2\bar{\xi}_0][(1-\xi_0\bar{\xi}_0)\xi_1-2\xi_0]}.
\]
Thus, using (\ref{e:key1}) and (\ref{e:key2}) we find
\begin{align}\nonumber
\frac{\eta_2\bar{\partial}_1\bar{\xi}_2+\bar{\eta}_2\bar{\partial}_1\xi_2}{(1+\xi_2\bar{\xi}_2)^2}=&\frac{\eta_1\bar{\partial}_1\bar{\xi}_1+\bar{\eta}_1\bar{\partial}_1\xi_1}{(1+\xi_1\bar{\xi}_1)^2}
+\frac{2(|\xi_0-\xi_1|^2-|1+\xi_0\bar{\xi}_1|^2)(\eta_0\bar{\partial}_1\bar{\xi}_0+\bar{\eta}_0\bar{\partial}_1\xi_0)}{(1+\xi_0\bar{\xi}_0)^3(1+\xi_1\bar{\xi}_1)}\\\nonumber
&+\frac{2[(\xi_0-\xi_1)(1+\xi_0\bar{\xi}_1)[(1+\xi_1\bar{\xi}_1)\bar{\partial}_1\bar{\xi}_0-(1+\xi_0\bar{\xi}_0)\bar{\partial}_1\bar{\xi}_1]]}{(1+\xi_0\bar{\xi}_0)^2(1+\xi_1\bar{\xi}_2)}\;r_0\\
&+\frac{2[(\bar{\xi}_0-\bar{\xi}_1)(1+\bar{\xi}_0\xi_1)[(1+\xi_1\bar{\xi}_1)\bar{\partial}_1\xi_0-(1+\xi_0\bar{\xi}_0)\bar{\partial}_1\xi_1]]}{(1+\xi_0\bar{\xi}_0)^2(1+\xi_1\bar{\xi}_2)}\;r_0\label{e:inta}.
\end{align}
Now, since we are reflecting in a surface, we have the potential
function $r_0$ satisfying
\[
\frac{2(\eta_0\bar{\partial}_0\bar{\xi}_0+\bar{\eta}_0\bar{\partial}_0\xi_0)}{(1+\xi_0\bar{\xi}_0)^2}=\bar{\partial}_0r_0,
\]
and so
\begin{align}
\frac{2(\eta_0\bar{\partial}_1\bar{\xi}_0+\bar{\eta}_0\bar{\partial}_1\xi_0)}{(1+\xi_0\bar{\xi}_0)^2}=&\frac{2(\eta_0(\bar{\partial}_1\bar{\nu}_0\bar{\partial}_0\bar{\xi}_0+\bar{\partial}_1\nu_0\partial_0\bar{\xi}_0)
+\bar{\eta}_0(\bar{\partial}_1\nu_0\partial_0\xi_0+\bar{\partial}_1\bar{\nu}_0\bar{\partial}_0\xi_0)}{(1+\xi_0\bar{\xi}_0)^2}\nonumber\\
&=\frac{2(\bar{\partial}_1\bar{\nu}_0(\eta_0\bar{\partial}_0\bar{\xi}_0+\bar{\eta}_0\bar{\partial}_0\xi_0)+\bar{\partial}_1\nu_0(\eta_0\partial_0\bar{\xi}_0+\bar{\eta}_0\partial_0\xi_0))}{(1+\xi_0\bar{\xi}_0)^2}\nonumber\\
&=\bar{\partial}_1\bar{\nu}_0\bar{\partial}_0r_0+\bar{\partial}_1\nu_0\partial_0r_0
\nonumber\\
&=\bar{\partial}_1r_0. \nonumber
\end{align}
Substituting this in (\ref{e:inta}) we get that
\begin{align}\nonumber
\frac{\eta_2\bar{\partial}_1\bar{\xi}_2+\bar{\eta}_2\bar{\partial}_1\xi_2}{(1+\xi_2\bar{\xi}_2)^2}=&\frac{\eta_1\bar{\partial}_1\bar{\xi}_1+\bar{\eta}_1\bar{\partial}_1\xi_1}{(1+\xi_1\bar{\xi}_1)^2}
+\frac{|\xi_0-\xi_1|^2-|1+\xi_0\bar{\xi}_1|^2}{(1+\xi_0\bar{\xi}_0)(1+\xi_1\bar{\xi}_1)}\bar{\partial}_1r_0\\\nonumber
&+\frac{2[(\xi_0-\xi_1)(1+\xi_0\bar{\xi}_1)[(1+\xi_1\bar{\xi}_1)\bar{\partial}_1\bar{\xi}_0-(1+\xi_0\bar{\xi}_0)\bar{\partial}_1\bar{\xi}_1]]}{(1+\xi_0\bar{\xi}_0)^2(1+\xi_1\bar{\xi}_2)}\;r_0\\\nonumber
&+\frac{2[(\bar{\xi}_0-\bar{\xi}_1)(1+\bar{\xi}_0\xi_1)[(1+\xi_1\bar{\xi}_1)\bar{\partial}_1\xi_0-(1+\xi_0\bar{\xi}_0)\bar{\partial}_1\xi_1]]}{(1+\xi_0\bar{\xi}_0)^2(1+\xi_1\bar{\xi}_2)}\;r_0,\nonumber
\end{align}
and the last three terms on the right hand side are $\bar{\partial}_1$ of
the real function, as stated in the theorem.
\end{pf}

As a corollary we get the Theorem of Malus: 

\begin{Cor}
A reflected congruence is integrable if and only if the initial
congruence is integrable.
\end{Cor}
\begin{pf}
Since the second term on the right-hand side of (\ref{e:twist}) is
$\bar{\partial}_1$ of a real function, the integrability condition
(\ref{e:intcond}) means that the outgoing congruence is integrable if
and only if the incoming wave is integrable.
\end{pf}

In the next sections we consider three special cases: that of
reflection of an arbitrary wavefront in a plane, and
that of plane and spherical wavefronts reflected off an arbitrary
surface. In all cases the
resulting wavefront can be explicitly determined. 

More generally, assume that the incident wave can be
described by $\eta_1=F_1(\xi_1,\bar{\xi}_1)$ for some complex function
$F_1$. Suppose further that we can invert (\ref{e:key1}) for $\xi_1$
as a function of $\xi_0$ and $\eta_0$. Then we can substitute
this in (\ref{e:reflawa}) and (\ref{e:key2})
to find both $\xi_2$ and $\eta_2$ as a function of $\xi_0$ and
$\eta_0$. Finally, we can find the potential function $r_2$ by
integrating (\ref{e:pot}).

\section{Reflection in a Plane}

Consider an arbitrary wave described parametrically by
$\xi_1(\nu,\bar{\nu})$ and $\eta_1(\nu,\bar{\nu})$. We want to determine
the resulting  wave  after reflection in a plane. By a rotation we can
align the plane with the $x^1x^2$-plane so that the initial wave lies in the region $x^3>0$.
\begin{Prop}\label{t:planesurf}
The reflection of a wave given by $\xi_1(\nu,\bar{\nu})$ and
$\eta_1(\nu,\bar{\nu})$ in the $x^1x^2$-plane is 
\begin{equation}\label{e:refinpl}
\xi_2=\frac{1}{\bar{\xi}_1},
\qquad\qquad
\eta_2=-\frac{\bar{\eta}_1}{\bar{\xi}_1^2}.
\end{equation}
\end{Prop}

\begin{pf}
The $x^1x^2$-plane is given by $\xi_0=r_0=0$. Inserting this in the
reflection law (\ref{e:reflawa}) immediately gives the first of the
above equations. Inserting it in (\ref{e:key1}) and (\ref{e:key2}) we
get
\begin{equation}\label{e:step1}
\eta_1=\eta_0-\xi_1^2\bar{\eta}_0,
\qquad\qquad
\eta_2=\eta_0-\frac{\bar{\eta}_0}{\bar{\xi}_1^2}.
\end{equation}
Adding the first of these to $\xi_1^2$ times its complex conjugate gives
\[
\eta_0=\frac{\eta_1+\xi_1^2\bar{\eta}_1}{1-\xi_1^2\bar{\xi}_1^2}.
\]
Inserting this in the second of (\ref{e:step1}) gives the result.
\end{pf}

\subsection{Reflection of a Spherical Wavefront in a Plane}
Consider a spherical wavefront with source ($0,0,t_1$). This
congruence is given by $\eta_1=-2t_1\xi_1$. Substituting this in 
(\ref{e:refinpl}) we find
\[
\xi_2=\frac{1}{\bar{\xi}_1}
\qquad\qquad
\eta_2=\frac{2t_1}{\bar{\xi}_1}.
\]
Combining these two we see that $\eta_2=2t_1\xi_2$, which is a
spherical congruence with source ($0,0,-t_1$). Thus we retrieve
the well-known law of reflection, whereby the reflection of a
spherical wave centered at ($0,0,t_1$) in a plane mirror is another
spherical wave with virtual centre ($0,0,-t_1$).

\section{Reflection of a Plane Wavefront}

Consider an incoming wavefront with fixed direction, i.e. a plane
wavefront. In this case, the reflected wave can be expressed
explicitly in terms of the reflecting surface:

\begin{Prop}\label{t:planepar}
The reflection of a plane wavefront with direction $\xi_1$ off a
surface given by $\xi_0(\nu,\bar{\nu})$ and $\eta_0(\nu,\bar{\nu})$ is
    given by
\begin{equation}\label{e:planewave1}
\xi_2=\frac{2\xi_0\bar{\xi}_1+1-\xi_0\bar{\xi}_0}
           {(1-\xi_0\bar{\xi}_0)\bar{\xi}_1-2\bar{\xi}_0},
\end{equation}
\begin{equation}\label{e:planewave2}
\eta_2=\frac{(\bar{\xi}_0-\bar{\xi}_1)^2\eta_0
-(1+\xi_0\bar{\xi}_1)^2\bar{\eta}_0
+(\bar{\xi}_0-\bar{\xi}_1)(1+\xi_0\bar{\xi}_1)(1+\xi_0\bar{\xi}_0)r_0
}{((1-\xi_0\bar{\xi}_0)\bar{\xi}_1-2\bar{\xi}_0)^2},
\end{equation}
the reflected congruence being parameterized by the parameter value
$\nu$ of the surface at the point of reflection.
\end{Prop} 
\begin{pf}

Since the incident angle is constant, there is no need to determine
the points of intersection, and the above equations are just the
equation of reflection (\ref{e:reflawa}) and equation (\ref{e:key2}). 
\end{pf}

In this proposition, the resulting wavefront is parameterized by the point
of reflection on the surface. For a plane wave it is also possible to parameterise the
outgoing wave by its direction:

\begin{Prop}\label{t:planedir}
Suppose $S$ a surface in ${\Bbb{R}}^3$, given by the
$\eta_0=F_0(\xi,\bar{\xi})$ with potential function
$r_0(\xi,\bar{\xi})$, is struck by a plane wave with normal direction
$-\bar{\xi}_1^{-1}$. Then the reflected wave is given by the line
congruence ($\xi$, $\eta_2=F_2(\xi,\bar{\xi})$) with
\begin{equation}\label{e:planeperp}
F_2(\xi,\bar{\xi})=\frac{1}{4}\left(
(1+\gamma)^2F_0(\xi_0,\bar{\xi}_0)
  -(\xi-\gamma\xi_1)^2\bar{F}_0(\xi_0,\bar{\xi}_0)
   +2(\xi_1-\xi)\gamma\;r_0(\xi_0,\bar{\xi}_0) 
\right),
\end{equation}
where
\begin{equation}\label{e:xi0}
\xi_0=\frac{\xi\bar{\xi}\xi_1\bar{\xi}_1-1\pm\sqrt{(1+\xi\bar{\xi})(1+\xi_1\bar{\xi}_1)(1+\bar{\xi}\xi_1)(1+\xi\bar{\xi}_1)}}
   {\bar{\xi}(1+\xi\bar{\xi}_1)+\bar{\xi}_1(1+\bar{\xi}\xi_1)},
\end{equation}
and
\[
\gamma=\pm\sqrt{\frac{(1+\xi\bar{\xi}_1)(1+\xi\bar{\xi})}
    {(1+\bar{\xi}\xi_1)(1+\xi_1\bar{\xi}_1)}}.
\]
The two signs yield the same line congruence with opposite
orientation, and can be chosen to point in the outgoing direction.
\end{Prop}
\begin{pf}
From the reflection law (\ref{e:reflawa}) with $\xi_2=\xi$ and
$\xi_1\rightarrow-\bar{\xi}_1^{-1}$ we have
\[
\xi=\frac{(\xi_0\bar{\xi}_0-1)\xi_1+2\xi_0}
           {2\bar{\xi}_0\xi_1+1-\xi_0\bar{\xi}_0}.
\]
This can be recast as a quadratic equation for $\xi_0$, which has solution
(\ref{e:xi0}). From these, after some calculation, we find the following expressions
\[
1+\xi_0\bar{\xi}_0=\pm\frac{2\sqrt{b^2+\beta\bar{\beta}}(-b\pm\sqrt{b^2+\beta\bar{\beta}})}{\beta\bar{\beta}},
\qquad\qquad
1-\xi_0\bar{\xi}_0=\frac{2b(-b\pm\sqrt{b^2+\beta\bar{\beta}})}{\beta\bar{\beta}},
\]
\[
1+\xi\bar{\xi}_0=\frac{\beta-b\xi\pm\xi\sqrt{b^2+\beta\bar{\beta}}}{\beta},
\qquad\qquad
\xi_0-\xi=\frac{-\xi\bar{\beta}-b\pm\sqrt{b^2+\beta\bar{\beta}}}{\bar{\beta}},
\]
where we have introduced
\[
b=1-\xi\bar{\xi}\xi_1\bar{\xi}_1 \qquad\qquad 
 \beta=\xi(1+\bar{\xi}\xi_1)+\xi_1(1+\xi\bar{\xi}_1).
\]
Then we compute the following combinations of the above
\[
\frac{1+\xi\bar{\xi}_0}{1+\xi_0\bar{\xi}_0}=
\frac{\pm\sqrt{(1+\xi\bar{\xi}_1)(1+\xi\bar{\xi})}+
        \sqrt{(1+\bar{\xi}\xi_1)(1+\xi_1\bar{\xi}_1)}}
 {2\sqrt{(1+\bar{\xi}\xi_1)(1+\xi_1\bar{\xi}_1)}},
\]
\[
\frac{\xi_0-\xi}{1+\xi_0\bar{\xi}_0}=
\frac{\pm\xi_1\sqrt{(1+\xi\bar{\xi}_1)(1+\xi\bar{\xi})}-
        \xi\sqrt{(1+\bar{\xi}\xi_1)(1+\xi_1\bar{\xi}_1)}}
 {2\sqrt{(1+\bar{\xi}\xi_1)(1+\xi_1\bar{\xi}_1)}},
\]
\[
\frac{(\xi_0-\xi)(1+\xi\bar{\xi}_0)}{1+\xi_0\bar{\xi}_0}=
\pm\frac{1}{2}(\xi_1-\xi)\sqrt{\frac{(1+\xi\bar{\xi}_1)(1+\xi\bar{\xi})}
   {(1+\bar{\xi}\xi_1)(1+\xi_1\bar{\xi}_1)}}.
\]
Substituting these in (\ref{e:key2}) gives the result
stated. 
\end{pf}

We also have the following corollary for waves travelling down the
$x^3$-axis:    

\begin{Cor}\label{c:planex3}
Suppose $S$ a surface in ${\Bbb{R}}^3$, given by the
$\eta_0=F_0(\xi,\bar{\xi})$ with potential function
$r_0(\xi,\bar{\xi})$, is struck by a plane wave moving down the
$x^3$-axis. Then the reflected wave is given by the line congruence
($\xi$, $\eta_2=F_2(\xi,\bar{\xi})$) with
\begin{equation}\label{e:planeperp1}
F_2(\xi,\bar{\xi})=\frac{1}{4}\left(
\left(1+\sqrt{1+\xi\bar{\xi}}\right)^2F_0(\xi_0,\bar{\xi}_0)
  -\xi^2\bar{F}_0(\xi_0,\bar{\xi}_0)
   -2\xi\sqrt{1+\xi\bar{\xi}}\;r_0(\xi_0,\bar{\xi}_0) 
\right),
\end{equation}
where
\[
\xi_0=\frac{-1+\sqrt{1+\xi\bar{\xi}}}{\bar{\xi}}.
\]
\end{Cor}
\begin{pf}
This follows from the above Proposition \ref{t:planedir} by setting $\xi_1=0$.
\end{pf}

By way of example, we now compute the twistor function for 
a plane wave  after reflection off the unit sphere and the rotationally symmetric torus.

\subsection{Reflection of a Plane Wavefront in the Unit Sphere}

By translation we can move the unit sphere to the origin in
${\Bbb{R}}^3$ and by a rotation we can set the direction of the plane
wave to be down the $x^3$-axis. Thus we can use corollary
\ref{c:planex3}. 

The twistor function and potential for the unit sphere at the origin
are $F_0=0$ and $r_0=1$. Thus, substituting these in
(\ref{e:planeperp1}), the reflected wavefront is given by 
\[
F_2=-\frac{1}{2}\xi\sqrt{1+\xi\bar{\xi}},
\qquad\qquad
r_2=\frac{2}{\sqrt{1+\xi\bar{\xi}}}+C.
\]
Figure 2 shows the resulting outgoing wavefronts (which are
axially symmetric, but not spherically symmetric).

\setlength{\epsfxsize}{3in}
\begin{center}
   \mbox{\epsfbox{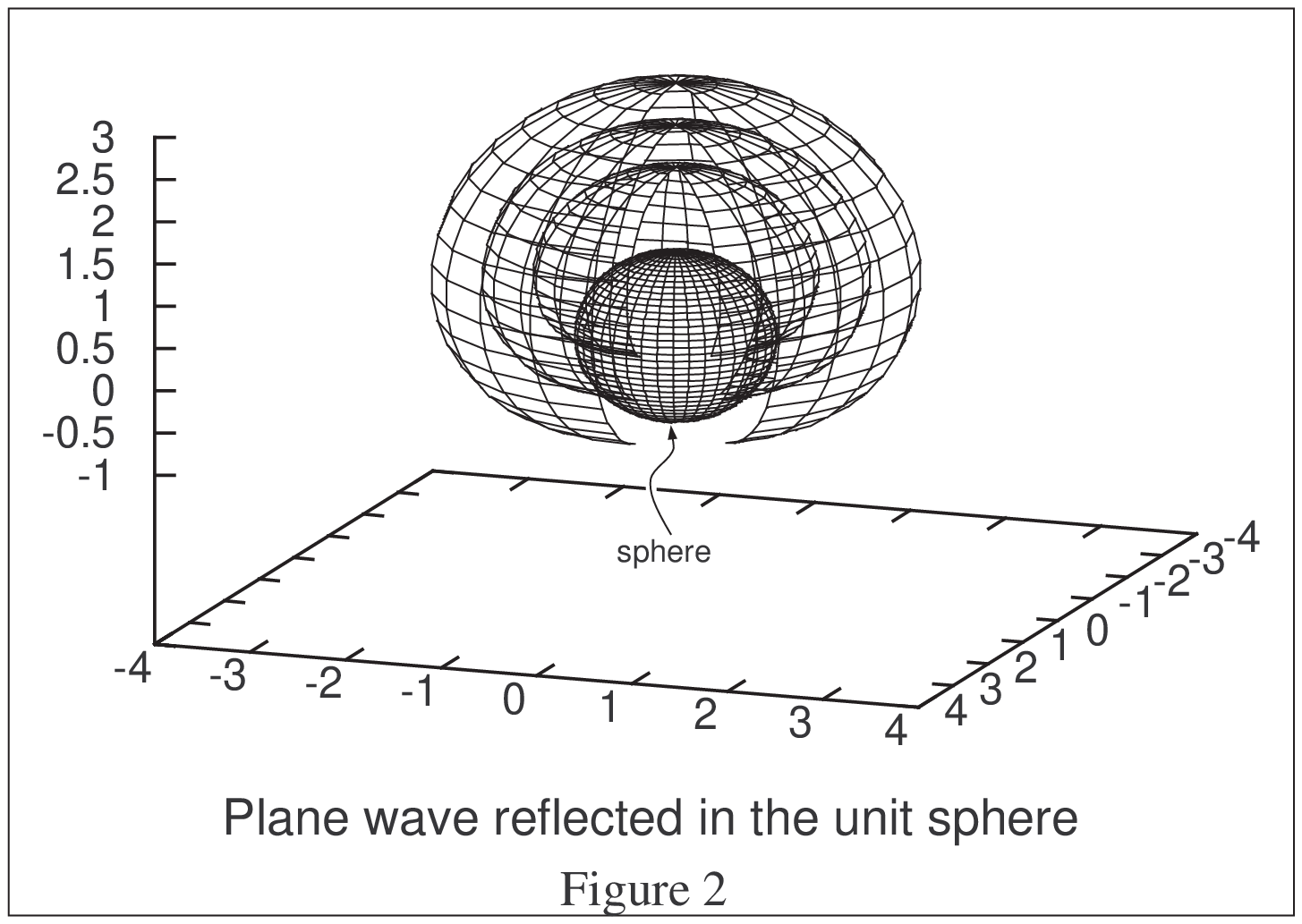}}
\end{center}

\subsection{Reflection of a Plane Wavefront in a Torus}

Consider the torus with core radius $a$ and meridian radius $b$ which is
rotationally symmetric about the $x^3$-axis. This has twistor function and potential
\[
F_0=\frac{a}{2}(1-\xi\overline{\xi})\sqrt{\frac{\xi}{\overline{\xi}}},
\qquad\qquad
r_0=\frac{2a\sqrt{\xi\overline{\xi}}}{1+\xi\overline{\xi}}+b.
\]
In what follows we will consider the case where $a=2$ and $b=1$. We
will deal with three cases: a wave travelling down the $x^3$-axis, up
the $x^1$-axis and one striking the torus at an angle of 45$^0$ to the
vertical. 

\subsubsection{Plane wavefront moving down the $x^3$-axis}
The line congruence resulting from reflection of a plane wave moving
down the $x^3$-axis can be found directly from corollary
\ref{c:planex3}. The result is
\[
F_2=\frac{2(1-\xi\overline{\xi})
             -\sqrt{\xi\bar{\xi}}\sqrt{1+\xi\bar{\xi}}}{2}\sqrt{\frac{\xi}{\overline{\xi}}},
\qquad\qquad
r_2=\frac{4\sqrt{\xi\overline{\xi}}+2\sqrt{1+\xi\bar{\xi}}}{1+\xi\overline{\xi}},
\] 
where we have integrated $F_2$ to obtain $r_2$.

Figure 3 shows the complete wavefront as it leaves the
torus. Notice the annular shadow.

\setlength{\epsfxsize}{3in}
\begin{center}
   \mbox{\epsfbox{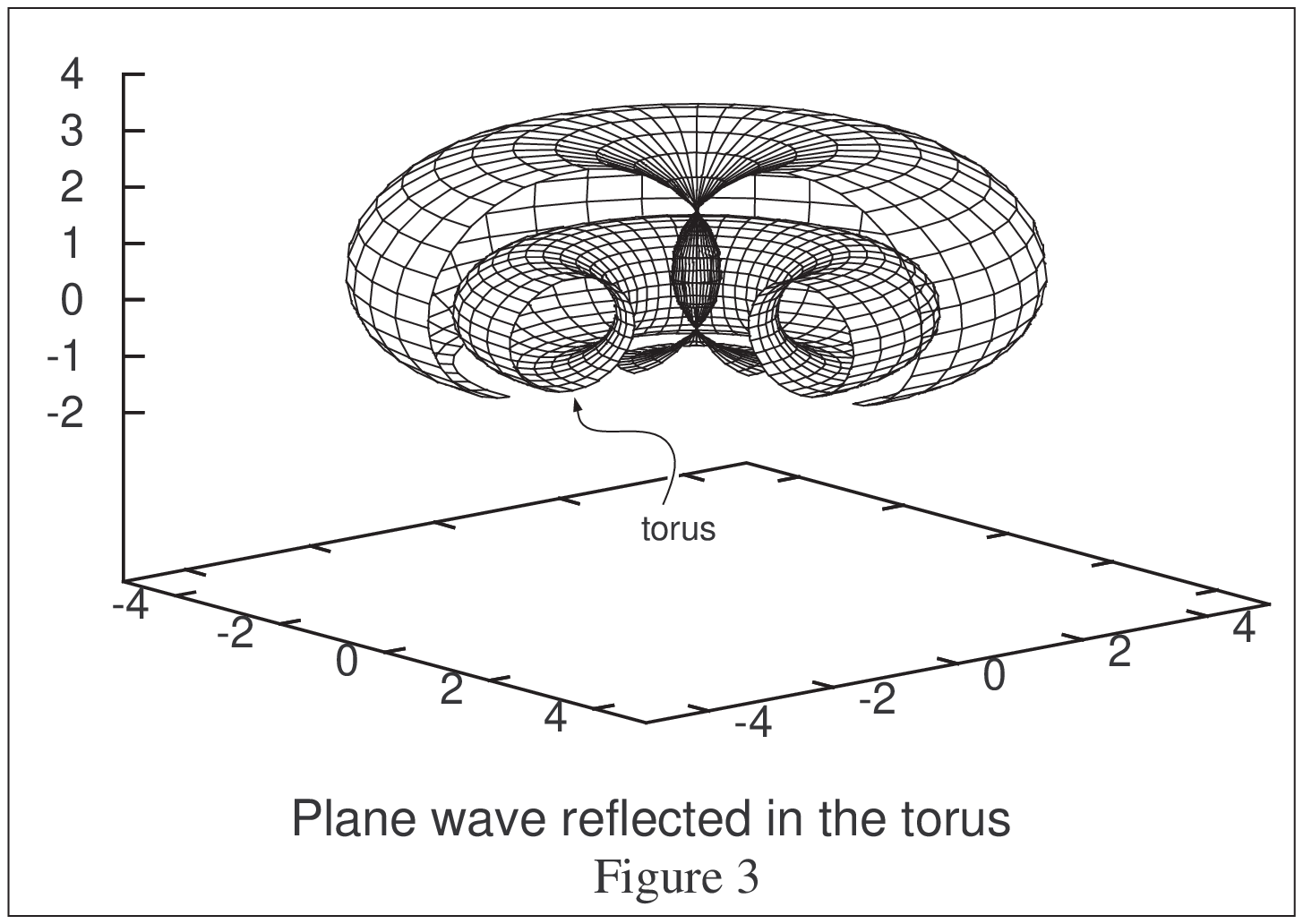}}
\end{center}

\subsubsection{Plane wavefront moving up the $x^1$-axis}
In order to determine the reflection of a plane wave approaching from
another direction we shall use Proposition \ref{t:planepar}. If the
direction of the incoming wave is $\xi_1$, substitution of the torus
twistor function and potential in (\ref{e:planewave2}) gives
\begin{equation}\label{e:pltorus1}
\eta_2=\frac{[(\bar{\xi}_0-\bar{\xi}_1)^2\xi_0
-(1+\xi_0\bar{\xi}_1)^2\bar{\xi}_0](1-\xi_0\bar{\xi}_0)
+(\bar{\xi}_0-\bar{\xi}_1)(1+\xi_0\bar{\xi}_1)
     [(1+\xi_0\bar{\xi}_0)\sqrt{\xi_0\bar{\xi}_0}+4\xi_0\bar{\xi}_0]}
{\sqrt{\xi_0\bar{\xi}_0}((1-\xi_0\bar{\xi}_0)\bar{\xi}_1-2\bar{\xi}_0)^2}.
\end{equation}
Integrating this gives the potential function
\begin{equation}\label{e:pltorus2}
r_2=\frac{4[2(|\bar{\xi}_0-\bar{\xi}_1|^2
-|1+\xi_0\bar{\xi}_1|^2)\sqrt{\xi_0\bar{\xi}_0}+
+(\xi_0\bar{\xi}_0(1-\xi_1\bar{\xi}_1)-\bar{\xi}_0\xi_1-\xi_0\bar{\xi}_1)
     (1+\xi_0\bar{\xi}_0)]}
{(1+\xi_1\bar{\xi}_1)(1+\xi_0\bar{\xi}_0)^2}.
\end{equation}

For a wave travelling along the $x^1$-axis we set $\xi_1=1$. The
resulting wavefront is shown in Figure 4.

\setlength{\epsfxsize}{3in}
\begin{center}
   \mbox{\epsfbox{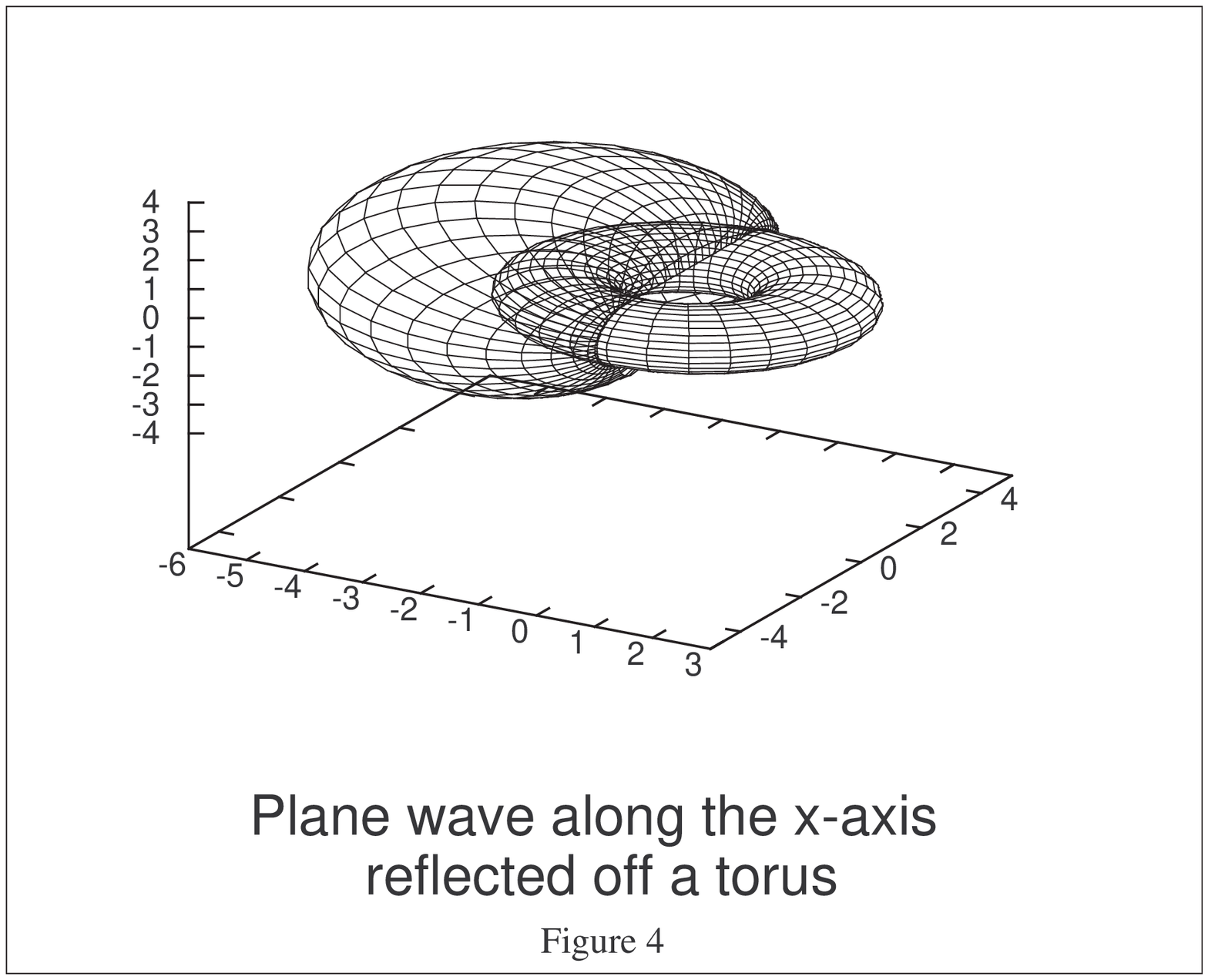}}
\end{center}

\subsubsection{Plane wavefront at an angle of 45$^0$ to the vertical}
Finally, a portion of the reflected wavefront generated by an incident
wave making an angle of 45$^0$ with the $x^3$-axis is shown in Figure 5. This is obtained by setting $\xi_1=2.4$ in (\ref{e:pltorus1})
and (\ref{e:pltorus2}).

\setlength{\epsfxsize}{3in}
\begin{center}
   \mbox{\epsfbox{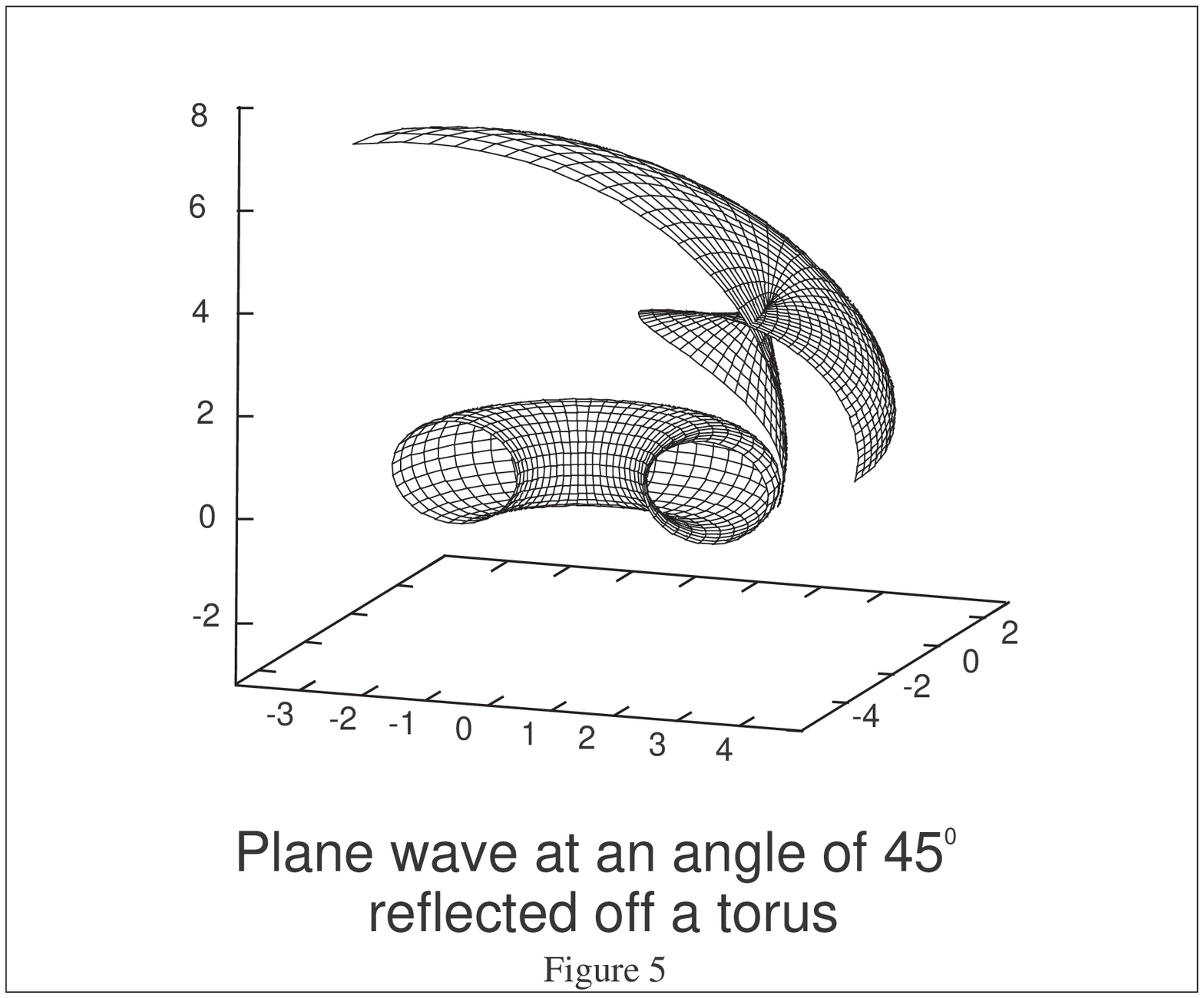}}
\end{center}

\section{Reflection of a Spherical Wavefront}

Assume that the wavefront is spherical, that is, the rays have a
single focus. By a translation we can move this focus to the origin
and then the wavefront is simply given by $\eta_1=0$. Every point in
${\Bbb{R}}^3$, other than the origin, is contained on two oriented
lines in the congruence (the same line with both orientations).

\begin{Prop}\label{t:sphwave}
A spherical wavefront, with focus the origin, reflected off a surface
(not containing the origin) determined by $\xi_0(\nu,\bar{\nu})$, $\eta_0(\nu,\bar{\nu})$ and $r_0(\nu,\bar{\nu})$ gives rise to the congruence
\begin{equation}\label{e:spherical1}
\xi_2=\frac{2\eta_0+2\bar{\eta}_0\xi_0^2\pm2\xi_0\beta_0}{2 (\xi_0\bar{\eta}_0-\bar{\xi}_0\eta_0)-(1+\xi_0\bar{\xi}_0)^2r_0\pm(1-\xi_0\bar{\xi}_0)\beta_0 },
\end{equation}
\begin{equation}\label{e:spherical2}
\eta_2=\alpha_1^2\eta_0-\alpha_2^2\bar{\eta}_0-2(1+\xi_0\bar{\xi}_0)\alpha_3r_0,
\end{equation}
where
\[
\alpha_1=\frac{2\bar{\eta}_0\xi_0-(1+\xi_0\bar{\xi}_0)r_0\pm\beta_0 }{2(\xi_0\bar{\eta}_0-\bar{\xi}_0\eta_0)
 -(1+\xi_0\bar{\xi}_0)^2r_0\pm(1-\xi_0\bar{\xi}_0)\beta_0 },
\]
\[
\alpha_2=\frac{2\eta_0+\xi_0(1+\xi_0\bar{\xi}_0)r_0\pm\xi_0\beta_0 }{2(\xi_0\bar{\eta}_0-\bar{\xi}_0\eta_0)
 -(1+\xi_0\bar{\xi}_0)^2r_0\pm(1-\xi_0\bar{\xi}_0)\beta_0 },
\]
\[
\alpha_3=\frac{4\eta_0\bar{\eta}_0\xi_0-(\eta_0-\bar{\eta}_0\xi_0^2)(1+\xi_0\bar{\xi}_0)r_0\pm(\eta_0+\bar{\eta}_0\xi_0^2)\beta_0 }{(2(\xi_0\bar{\eta}_0-\bar{\xi}_0\eta_0)
 -(1+\xi_0\bar{\xi}_0)^2r_0\pm(1-\xi_0\bar{\xi}_0)\beta_0 )^2},
\]
and
\[
\beta_0=\sqrt {4 \eta_0 \bar{\eta}_0+(1+\xi_0\bar{\xi}_0)^2r_0^2}.
\]
Here the two different signs give the same line congruence with the
opposite orientation, and can be chosen to point in the outgoing direction.

\end{Prop}
\begin{pf}

Setting $\eta_1=0$ in (\ref{e:key1}) we get a quadratic equation for
$\xi_1$, with solution:
\[
\xi_1=\frac{2(\eta_0\bar{\xi}_0+\bar{\eta}_0\xi_0)
 -(1-\xi^2_0\bar{\xi}_0^2)r_0 \pm(1+\xi_0\bar{\xi}_0)\beta_0}
{2(\bar{\eta}_0-\eta_0\bar{\xi}_0^2)+2\bar{\xi}_0(1+\xi_0\bar{\xi}_0)r_0}.
\]  
Inserting this in (\ref{e:reflawa}) yields the expression for $\xi_2$
in the theorem. In addition we find that
\[
\bar{\xi}_1-\bar{\xi}_0=\frac{(2\bar{\eta}_0\xi_0-(1+\xi_0\bar{\xi}_0)r_0\pm\beta_0 )(1+\xi_0\bar{\xi}_0)}{2(\eta_0-\bar{\eta}_0\xi_0^2)+2\xi_0(1+\xi_0\bar{\xi}_0)r_0},
\]
\[
1+\xi_0\bar{\xi}_1=\frac{(2\eta_0+\xi_0(1+\xi_0\bar{\xi}_0)r_0\pm\xi_0\beta_0 )(1+\xi_0\bar{\xi}_0)}{2(\eta_0-\bar{\eta}_0\xi_0^2)+2\xi_0(1+\xi_0\bar{\xi}_0)r_0},
\]
while
\[
(1-\xi_0\bar{\xi}_0)\bar{\xi}_1-2\bar{\xi}_0=\frac{(2(\xi_0\bar{\eta}_0-\bar{\xi}_0\eta_0)
 -(1+\xi_0\bar{\xi}_0)^2r_0\pm(1-\xi_0\bar{\xi}_0)\beta_0 )(1+\xi_0\bar{\xi}_0)}{2(\eta_0-\bar{\eta}_0\xi_0^2)+2\xi_0(1+\xi_0\bar{\xi}_0)r_0}.
\]
Substituting these in (\ref{e:key2}) completes the theorem.
\end{pf}

\subsection{Reflection of a Spherical Wavefront in the Unit Sphere}

Consider the unit sphere, centred at (0,0,$-$2). This is determined by

\[
\eta_0=2\xi_0, \qquad \qquad r_0=1-2\frac{1-\xi_0\bar{\xi}_0}{1+\xi_0\bar{\xi}_0}.
\]
Thus
\[
\beta_0=\sqrt{1+10\xi_0\bar{\xi}_0+r\xi_0^2\bar{\xi}_0^2},
\]
and by (\ref{e:spherical1}) 
\[
\xi_2=\frac{2\xi_0(2(1+\xi_0\bar{\xi}_0)+\beta_0)}
  {1-2\xi_0\bar{\xi}_0-3\xi_0^2\bar{\xi}_0^2+(1-\xi_0\bar{\xi}_0)\beta_0},
\]
where we have chosen the $+$ sign in the equations. Further
straightforward, if lengthy, calculations establish that
\[
\alpha_1=\frac{1+\xi_0\bar{\xi}_0+\beta_0}{1-2\xi_0\bar{\xi}_0-3\xi_0^2\bar{\xi}_0^2+\beta_0},
\qquad\qquad
\alpha_2=\frac{\xi_0(3(1+\xi_0\bar{\xi}_0)+\beta_0)}{1-2\xi_0\bar{\xi}_0-3\xi_0^2\bar{\xi}_0^2+\beta_0},
\]
\[
\alpha_3=\frac{\xi_0(1+\xi_0\bar{\xi}_0+\beta_0)}{(1-2\xi_0\bar{\xi}_0-3\xi_0^2\bar{\xi}_0^2+\beta_0)^2}.
\]
Finally, substitution of these in (\ref{e:spherical2}) gives
\[
\eta_2=\frac{4 \xi_0(1-3\xi_0\bar{\xi}_0)(1+3 \xi_0\bar{\xi}_0+ \sqrt {1+10 \xi_0\bar{\xi}_0+9 \xi_0^{2} \bar{\xi}_0^{2}} )}{1 +\xi_0 \bar{\xi}_0 -7 \xi_0^{2} \bar{\xi}_0^{2}+9 \xi_0^{3} \bar{\xi}_0^{3}+(1-4\xi_0 \bar{\xi}_0+3 \xi_0^{2} \bar{\xi}_0^{2})\sqrt {1+10 \xi_0 \bar{\xi}_0+ 9 \xi_0^{2} \bar{\xi}_0^{2}} }.
\]
To find the potential function $r_2$ we invert (\ref{e:pot}) with
$\nu=\xi_0$, since we are parameterising the outgoing rays by the
direction of the normal to the surface at the point of
intersection. After some computation we obtain
\[
r_2=\frac{-2(1-3\xi_0 \bar{\xi}_0)^2\sqrt {1+10 \xi_0 \bar{\xi}_0+ 9
    \xi_0^{2} \bar{\xi}_0^{2}}}{1+11 \xi_0 \bar{\xi}_0+19 \xi_0^{2} \bar{\xi}_0^{2}+ 9\xi_0^{3} \bar{\xi}_0^{3}}+C.
\]

At this juncture we have found parametric expressions for $\xi_2$,
$\eta_2$ and $r_2$, that is, the complete description of the outgoing
wave. If we graph the wavefronts the results are similar to the plane
wave reflected off the sphere. By way of comparison, Figure 6
compares a cross-section of the two wavefronts in a plane containing
the $x^3$-axis. The shaded area represents the shadow cast by the surface.

\setlength{\epsfxsize}{4in}
\begin{center}
   \mbox{\epsfbox{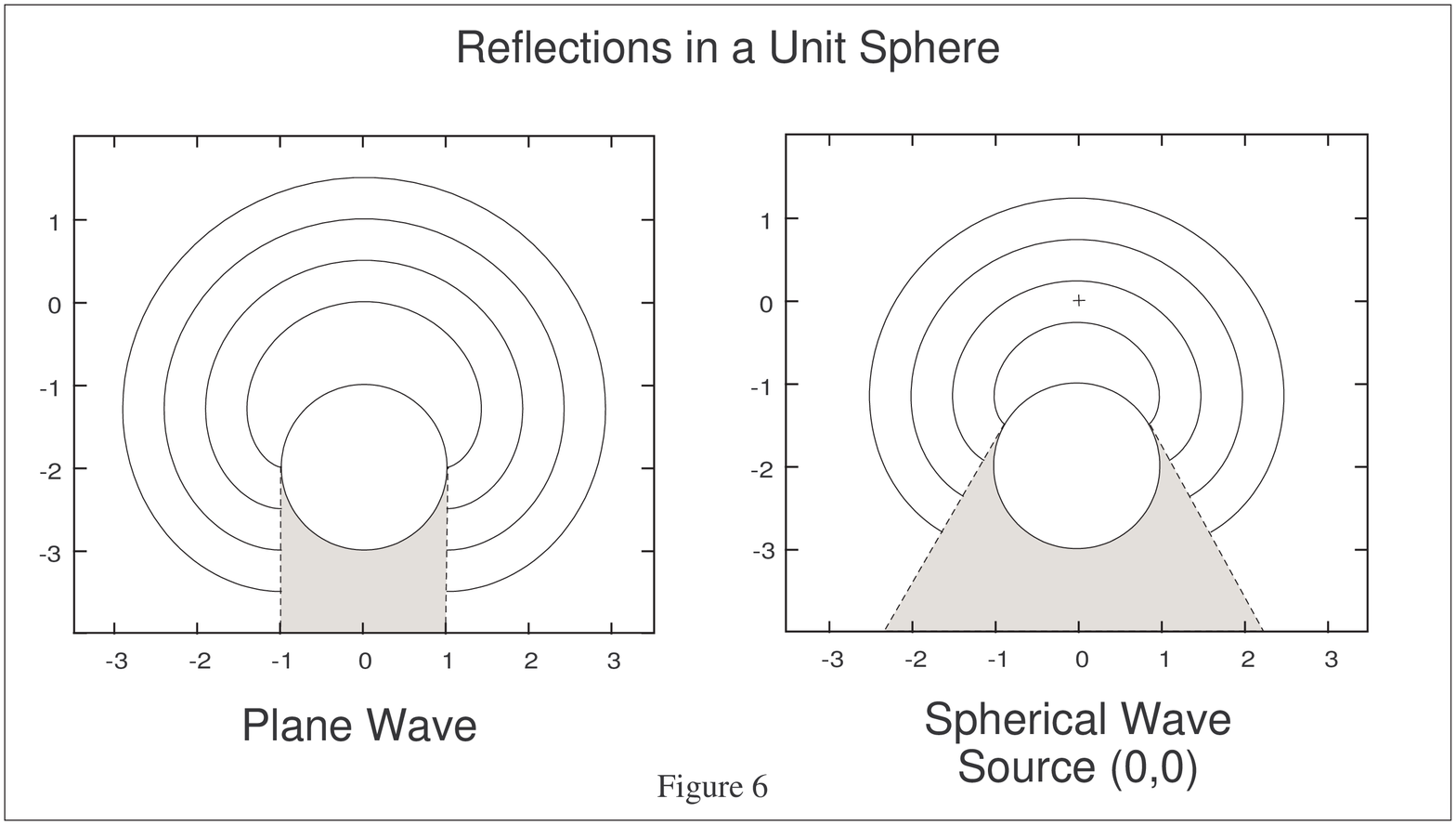}}
\end{center}

\end{document}